\documentclass[12pt, a4paper]{amsart}

\usepackage[colorlinks, linkcolor=blue, anchorcolor=blue, citecolor=green]{hyperref}

\usepackage{graphics,epic}
\usepackage{amsmath,amssymb, amsthm,mathrsfs}
\usepackage[all,2cell]{xy}
\usepackage{shorttoc}
\usepackage{tikz}

\newtheorem{theorem}{Theorem}[section]
\newtheorem*{theorem*}{Theorem}

\newtheorem{lemma}[theorem]{Lemma}

\newtheorem{corollary}[theorem]{Corollary}

\newtheorem*{conjecture*}{Conjecture}

\newtheorem{remark}[theorem]{Remark}
\newtheorem{definition}[theorem]{Definition}

\newtheorem{thm}[theorem]{Theorem}
\newtheorem{lem}[theorem]{Lemma}
\newtheorem{prop}[theorem]{Proposition}

\newcommand{\ie}{{\em i.e.}\ }

\newcommand{\opname}[1]{\operatorname{\mathsf{#1}}}

\renewcommand{\mod}{\opname{mod}\nolimits}

\newcommand{\add}{\opname{add}\nolimits}

\newcommand{\der}{\md}

\newcommand{\dimv}{\underline{\dim}\,}
\renewcommand{\Im}{\opname{Im}\nolimits}

\newcommand{\ind}{\opname{ind}\nolimits}

\newcommand{\Sub}{\opname{Sub}}
\newcommand{\Fac}{\opname{Fac}}

\newcommand{\D}{\mathbb{D}}

\renewcommand{\P}{\mathbb{P}}

\newcommand{\Hom}{\opname{Hom}}

\newcommand{\Ext}{\opname{Ext}}

\newcommand{\End}{\opname{End}}

\newcommand{\R}{{\mathbb R}}

\newcommand{\mc}{\mathcal{C}}
\newcommand{\md}{\mathcal{D}}
\newcommand{\mf}{\mathcal{F}}
\newcommand{\mh}{\mathcal{H}}

\newcommand{\mt}{\mathcal{T}}

\newcommand{\mx}{\mathcal{X}}
\newcommand{\my}{\mathcal{Y}}

\newcommand{\tauni}{\tau^{-1}}

\newcommand{\pd}{\operatorname{pd}\nolimits}
\newcommand{\id}{\operatorname{id}\nolimits}
\newcommand{\gld}{\operatorname{gldim}\nolimits}

\newcommand{\mfd}{\mathfrak{D}}

\newcommand{\mfr}{\mathfrak{R}}
\newcommand{\mfc}{\mathfrak{C}}

\begin{document}

\title[$\tau$-tilting graphs of quasi-silted algebras]{On $\tau$-tilting graphs for quasi-silted algebras}\thanks{Partially supported by the National Natural Science Foundation of China (Grant No. 11971326, 12171397, 12471037, 12571040)}

\author[Dai]{Wei Dai}
\address{Wei Dai\\ Department of Mathematics\\SiChuan University\\610064 Chengdu\\P.R.China}
\email{daiweimath@gmail.com}
\author[Fu]{Changjian Fu}
\address{Changjian Fu\\Department of Mathematics\\SiChuan University\\610064 Chengdu\\P.R.China}
\email{changjianfu@scu.edu.cn}
\author[Geng]{Shengfei Geng}
\address{Shengfei Geng\\Department of Mathematics\\SiChuan University\\610064 Chengdu\\P.R.China}
\email{genshengfei@scu.edu.cn}
\author[Liu]{Pin Liu}
\address{Pin Liu\\Department of Mathematics\\
   Southwest Jiaotong University\\
  610031 Chengdu \\
  P.R.China}
  \email{pinliu@swjtu.edu.cn}
 
\subjclass[2010]{}

\dedicatory{Dedicated to Professor Bangming Deng on the Occasion of his 60th Birthday}

\begin{abstract}
We prove that the $\tau$-tilting graph of any quasi-silted algebra is connected and has the reachable-in-face property. Our approach utilizes $\tau$-reduction and wall and chamber structures. In particular, we observe a sufficient condition on the wall and chamber structure under which the connectivity of $\tau$-tilting graphs is preserved under taking quotients of algebras. As an immediate consequence, the connectivity of $\tau$-tilting graphs is also established for several new classes of algebras.
\end{abstract}
\maketitle

\section{Introduction}
An important combinatorial invariant of a cluster algebra is its exchange graph, whose vertices correspond to seeds and whose edges connect seeds related by a single mutation. A cluster algebra is said to be of finite type if its exchange graph is finite, meaning it possesses only finitely many distinct seeds. Cluster algebras of finite type were classified in \cite{FZ03a} and shown to correspond to finite root systems. Moreover, the exchange graph of a cluster algebra of finite type can be realized as the 1-skeleton of a generalized associahedron, or Stasheff’s polytope \cite{FZ03b}. Through the representation-theoretic categorification of cluster algebras, one obtains a wide variety of exchange graphs associated with finite-dimensional algebras or differential graded (dg) algebras concentrated in non-positive degrees. These constructions stem from variations of tilting theory, where the vertices of the resulting exchange graphs correspond to support $\tau$-tilting modules, torsion pairs, silting objects, and related structures. For a comprehensive overview of the history of abstract exchange graphs arising from representation theory, we refer the reader to the introduction of \cite{BY2013}.

The exchange graph for a finite-dimensional algebra $A$ is also referred to as the $\tau$-tilting graph of $A$. In this graph, the vertices correspond to basic support $\tau$-tilting pairs, and the edges connect pairs related by a single mutation. It can be viewed as the dual graph of the abstract polytope $\Delta(A)$ of $A$ arising from $\tau$-tilting theory \cite{DIJ19}. The combinatorial properties of this abstract polytope have attracted significant attention. One of the fundamental questions in this area is to determine the number of connected components of the $\tau$-tilting graph for a given algebra.

It is well known that the $\tau$-tilting graph of $A$ is connected if $A$ belongs to one of the following classes of algebras:\begin{itemize}\item[(1)] algebras whose $\tau$-tilting graph has a finite connected component—in particular, algebras with only finitely many basic support $\tau$-tilting pairs \cite{AIR};\item[(2)] cluster-tilted algebras arising from hereditary abelian categories \cite{BMRRT,BKL,FG};\item[(3)] $2$-Calabi-Yau tilted algebras arising from marked surfaces, excluding closed surfaces with exactly one puncture \cite{QZ,Y};\item[(4)] gentle algebras \cite{FGLZ} and skew-gentle algebras \cite{HZZ22}.\end{itemize}

The main objective of this note is to show that the $\tau$-tilting graph of any quasi-silted algebra is connected (Theorem \ref{thm:connected of quasi-silted algebra}) and has the so-called reachable-in-face property (Corollary \ref{cor:reachable-in-face}), thereby adding a prominent new class of finite-dimensional algebras to the above list. Quasi-silted algebras were introduced in \cite{BZ16} as a generalization of quasi-tilted algebras \cite{HRS}, and they constitute an important class of finite-dimensional algebras. Notably, it was proven in \cite{BZ16} that quasi-silted algebras coincide with the shod algebras introduced in \cite{CL99}.

Our approach utilizes $\tau$-reduction \cite{Jasso} and the wall-and-chamber structure \cite{BST}. Along the way, we establish a sufficient condition under which this connectivity is preserved under taking quotients of algebras. As an immediate consequence, we obtain new connectivity results that significantly expand the known classes of algebras with connected $\tau$-tilting graphs, including, in particular, quotient algebras of skew-gentle algebras and quotient algebras of tame hereditary algebras.

\subsection*{Convention}
Throughout this paper, $k$ denotes an algebraically closed field. We write $\D = \Hom_k(-, k)$ for the $k$-duality. All algebras considered in this paper are assumed to be basic and finite-dimensional over $k$. By a hereditary abelian category, we mean a $k$-linear hereditary abelian category with finite-dimensional $\Hom$ and $\Ext$ spaces.

For a finite-dimensional $k$-algebra $A$, we denote by $\gld A$ its global dimension, and $\mod A$ the category of finitely generated right $A$-modules. For a module $M \in \mod A$, we denote by $|M|$ the number of pairwise non-isomorphic indecomposable direct summands of $M$. Furthermore, we denote by $\pd_A M$ (resp., $\id_A M$) the projective (resp., injective) dimension of $M$ in $\mod A$, and by $\opname{Fac} M$ the full subcategory of $\mod A$ consisting of all factor modules of finite direct sums of copies of $M$.

\section{Preliminary}
\subsection{(Support) $\tau$-tilting graphs}
 Let $A$ be a finite-dimensional $k$-algebra.  Let $\tau_A$, simply denoted by $\tau$, be the Auslander-Reiten translation of $\mod A$. 
 Recall that a module $M\in \mod A$ is {\it $\tau$-rigid} if $\Hom_A(M,\tau M)=0$. A $\tau$-rigid module $M$ is {\it $\tau$-tilting} if $|M|=|A|$.
A {\it $\tau$-rigid pair} is a pair $(M,P)$ with $M\in \mod A$ and $P$ a finitely generated projective right $A$-module, such that $M$ is $\tau$-rigid and $\Hom_A(P,M)=0$.
A $\tau$-rigid pair $(M,P)$ is a {\it  $\tau$-tilting pair} provided that $|M|+|P|=|A|$.
In this case, $M$ is a {\it support $\tau$-tilting} $A$-module and $P$ is uniquely determined by $M$ provided that $P$ is basic. In the following, we always identify basic support $\tau$-tilting modules with  basic $\tau$-tilting pairs. Denote by $\opname{\tau-tilt} A$ the set of all basic $\tau$-tilting pairs of $A$. 

Let $(M,P)$ and $(N,Q)$ be two $\tau$-rigid pairs. We say that $(N,Q)$ is a \emph{direct summand} of $(M,P)$ if $N$ and $Q$ are direct summands of $M$ and $P$, respectively. A $\tau$-rigid pair $(M,P)$ is said to be \emph{indecomposable} if $|M|+|P|=1$. In particular, each $\tau$-tilting pair has exactly $|A|$ pairwise non-isomorphic indecomposable direct summands. It is also known that every basic $\tau$-rigid pair is a direct summand of a  $\tau$-tilting pair.

Let $(M, P)$ be a basic $\tau$-rigid pair such that $|M|+|P|=|A|-1$. It was proved in \cite{AIR} that there exist exactly two non-isomorphic basic $\tau$-tilting pairs, $(M_1,P_1)$ and $(M_2,P_2)$, such that $(M,P)$ is a direct summand of $(M_i, P_i)$ for $i=1,2$ (cf.\ also \cite{DF15}). Clearly, $(M_1, P_1)$ and $(M_2, P_2)$ differ by exactly one indecomposable direct summand, say $(N,Q)$. In this case, $(M_1,P_1)$ is called the \emph{mutation} of $(M_2,P_2)$ at $(N,Q)$ and vice versa, which is denoted simply by $\mu_{(N,Q)}(M_1,P_1)=(M_2,P_2)$.

For two $\tau$-tilting pairs $(X,P)$ and $(Y,R)$, we write $(X,P)\sim (Y,R)$ if they are related by a sequence of mutations. 
\begin{definition}
The {\it $\tau$-tilting graph} $\mh(\opname{\tau-tilt} A)$ has vertex set indexed by the isomorphism classes of basic support $\tau$-tilting $A$-modules, and any two basic support $\tau$-tilting modules are connected by an edge if and only if they are mutations of each other.
\end{definition}

Demonet, Iyama and Jasso \cite{DIJ19} introduced a simplicial complex $\Delta(A)$ for $A$. In particular, there is a one-to-one correspondence between the $d$-simplexes of $\Delta(A)$ and the basic $\tau$-rigid pairs of $A$ which
have exactly $d+1$ indecomposable summands. The simplicial complex $\Delta(A)$ has pure
dimension of $|A|-1$. The  $\tau$-tilting graph $\mh(\opname{\tau-tilt} A)$ can be identified with the
dual graph of $\Delta(A)$. Each basic $\tau$-rigid pair $(L,R)$ determines
a face $\mathcal{H}_{(L,R)}(\opname{\tau-tilt} A)$ of $\mh(\opname{\tau-tilt} A)$. Namely, the face $\mathcal{H}_{(L,R)}(\opname{\tau-tilt} A)$ is the full subgraph of $\mh(\opname{\tau-tilt} A)$
consisting of basic support $\tau$-tilting pairs which admit $(L,R)$ as a direct summand. We also denote by
 $\tau$-tilt$_{(L,R)}A$ the set of $\tau$-tilting $A$-pairs containing $(L,R)$ as a direct summand.

\begin{definition}\cite{FGLZ}
The $\tau$-tilting graph $\mh(\opname{\tau-tilt} A)$ is said to have the 
{\em reachable-in-face} property if for any pair of basic $\tau$-tilting pairs 
$(M,P) \sim (N,Q)$ sharing a common direct summand $(L,R)$, 
there exists a path from $(M,P)$ to $(N,Q)$ that lies entirely within the face 
$\mathcal{H}_{(L,R)}(\opname{\tau-tilt} A)$ determined by $(L,R)$. In 
this case, we also say that the algebra $A$ has the reachable-in-face property. 
\end{definition}

It was observed in \cite{FGLZ} that there is a close connection between the reachable-in-face property and the connectedness of the $\tau$-tilting graph. While it is not yet known whether the reachable-in-face property holds for arbitrary $\tau$-tilting graphs, the following result provides an indication.

\begin{lem}\cite[Corollary 3.12]{CWZ}\label{l:cwz}
Let $(M,P)$ be a basic $\tau$-tilting pair of $A$ and $(Q, 0)$ a basic $\tau$-rigid pair which is a common direct summand of $(M,P)$ and $(A,0)$. If there is a path from $(A,0)$ to $(M,P)$ in $\mh(\opname{\tau-tilt} A)$, then there also exists such a path in the subgraph $\mathcal{H}_{(Q, 0)}(\opname{\tau-tilt} A)$.
\end{lem}

Let $\mathcal{T}$ be a functorially finite torsion class of $\mod A$. Recall that an object $X\in\mathcal{T}$ is $\Ext$-projective if $\Ext^1_A(X, \mathcal{T})=0$. Denote by $P(\mathcal{T})$ the direct sum of one copy of each of the indecomposable $\Ext$-projective objects in $\mathcal{T}$ up to isomorphism. In particular, $P(\mathcal{T})$ is a basic support $\tau$-tilting $A$-module. 
Conversely, for a basic $\tau$-rigid pair $(U,R)$, $^\perp\tau U\cap R^\perp$ is a functorially finite torsion class of $\mod A$ and $U$ is a direct summand of $P(^\perp\tau U\cap R^\perp)$ (cf. \cite{AIR}). 

Let $(U,R)$ be a basic $\tau$-rigid pair of $A$.
Denote by $\opname{\tau-rigid-pair} A$ the set of isomorphism classes of basic $\tau$-rigid pairs of $A$ and $\opname{\tau-rigid-pair}_{(U,R)}A$ the subset of $\opname{\tau-rigid-pair} A$ consisting of basic $\tau$-rigid pairs which admit $(U,R)$ as a direct summand. Let $T=P(^\perp\tau U\cap R^{\perp})$ and $A_{(U,R)}:=\End (T)/\langle e_U \rangle$, where $e_U$ is the idempotent of $\End(T)$ associated with $U$. The algebra $A_{(U,R)}$ is called the {\em $\tau$-reduction} of $A$ with respect to $(U,R)$.

The following is known as the $\tau$-reduction theory of $A$ with respect to $(U,R)$ (cf. \cite[Theorem 3.16]{Jasso} and \cite[Corollary A.4]{FGLZ}).
\begin{lemma}\label{l:tau-reduction}
    Keep the notation as above. There is an order-preserving bijection \[E_{(U,R)}: \opname{\tau-rigid-pair}_{(U,R)}A\to \opname{\tau-rigid-pair} A_{(U,R)}\] which commutes with direct sums and restricts to a bijection $$E_{(U,R)}:\opname{\tau-tilt}_{(U,R)}A\to \opname{\tau-tilt} A_{(U,R)}$$ commuting with mutations. In particular, there is an isomorphism of graphs between $\mathcal{H}_{(U,R)}(\opname{\tau-tilt}A)$ and $\mathcal{H}(\opname{\tau-tilt}A_{(U,R)})$.
\end{lemma}

\subsection{Wall and chamber structure}
We now recall a construction from \cite{BST}.
Let $A$ be a finite-dimensional $k$-algebra and $\{e_1,\dots, e_n\}$ a complete set of pairwise orthogonal idempotents of $A$. Let $P(i)=e_iA$ be the indecomposable projective $A$-module associated with $e_i$ and $S_i=\opname{top} P(i)$ its simple top, where $1\leq i\leq n$. 
 We identify the Grothendieck group $K_0(\mod A)$ of $\mod A$  with $\mathbb{Z}^n$ via the function $$\dimv: \mod A\to \mathbb{Z}^n$$ which maps $S_i$ to $\mathbf{e}_i$, where $\mathbf{e}_1,\dots, \mathbf{e}_n$ is the standard basis of $\mathbb{Z}^n$. Denote by $\langle-,-\rangle$ the canonical inner product of $\mathbb{R}^n$.
 For any vector $\theta\in \R^n$, a nonzero $A$-module $M$ is called {\it $\theta$-semistable} if $\langle \theta, \dimv M\rangle =0$ and $\langle \theta, \dimv L\rangle\leq 0$ for every submodule $L$ of $M$. The {\it stability space} of an $A$-module $M$ is then defined as
\[\mfd_A(M)=\{\theta\in \R^n~|~\text{$M$ is $\theta$-semistable}\}.\]
We say that $\mfd_A(M)$ is a {\it wall} of $A$ when $\mfd_A(M)$ has codimension one.

Outside the walls, there are only vectors $\theta$ having no nonzero $\theta$-semistable modules.
Removing the closure of all walls, we obtain a set
$$\mfr_A=\R^n\backslash \overline{\bigcup_{M\in\mod A}\mfd_A(M)}$$
whose connected components $\mfc$ are called {\it chambers}. As connected components of an
open set in $\R^n$, the chambers have dimension $n$. This decomposition of $\R^n$ is called the
{\it wall and chamber structure} of the algebra $A$ on $\R^n$.

 The following is an easy observation (cf. \cite[Lemma 4.13]{BST}).
\begin{lem}\label{l:B wall is A wall}
Let $A$ be a finite-dimensional $k$-algebra and $B$ a quotient algebra of $A$ with $|B|=|A|$.  Then every wall of $B$ is also a wall of $A$.
 \end{lem}

Let $(M,P)$ be a $\tau$-rigid pair. We decompose $P$ into indecomposable projective $A$-modules as: 
\[
P=\bigoplus_{i=1}^nP(i)^{c_i},
\]where $c_1,\dots, c_n$ are non-negative integers. Let $$\bigoplus_{i=1}^nP(i)^{b_i}\xrightarrow{f_M}\bigoplus_{i=1}^nP(i)^{a_i}\to M\to 0$$ be a minimal projective presentation of $M$, where $a_i,b_j$ are non-negative integers. Recall that the {\it $g$-vector} $g_{(M,P)}$ associated with $(M,P)$ is defined as
 \[
 g_{(M,P)}=[a_1-b_1,\dots, a_n-b_n]^t-\sum_{i=1}^nc_i\mathbf{e}_i.
 \]
 It is known that different $\tau$-rigid pairs have different $g$-vectors and the $n$ $g$-vectors of indecomposable direct summands of a basic $\tau$-tilting pair form a basis of $\mathbb{Z}^n$ (cf. \cite{AIR}).

The following fact was first noticed in \cite{BST} and was later shown in \cite{Asai}.
\begin{thm}\label{l:mutation and neighboring  chamber}
Let $A$ be a finite-dimensional $k$-algebra. Then there is an injective function $\mfc$ mapping the basic $\tau$-tilting pair $(M,P)$ onto a chamber $\mfc_{(M,P)}$ and every chamber arises this way. Moreover, $\tau$-tilting pair $(M',P')$ is a mutation of $(M,P)$ if and only if $\mfc_{(M',P')}$ is a neighbor of $\mfc_{(M,P)}$, namely, they are separated by a wall.  
\end{thm}
More precisely, let $(M,P)=\bigoplus_{i=1}^n(M_i,P_i)$ be a basic $\tau$-tilting pair with indecomposable direct summands $(M_i,P_i), 1\leq i\leq n$. The chamber $\mfc_{(M,P)}$ associated with $(M,P)$ is defined as
 \[
 \mfc_{(M,P)}=\left\{\sum_{i=1}^nk_ig_{(M_i,P_i)}~|~0<k_i\in \mathbb{R}\right\}.
 \]
 In other words, $\mfc_{(M,P)}$ is the interior of the positive cone \[{C_{(M,P)}}=\left\{\sum_{i=1}^nk_ig_{(M_i,P_i)}~|~0\leq k_i\in \mathbb{R}\right\}.\]

Recall that a smooth path $\gamma:[0,1]\rightarrow \R^n$ is a
{\it $\mfd_A$-generic path} if:
\begin{itemize}
\item  $\gamma(0)$ and $\gamma(1)$  are located inside some chambers;
\item  If $\gamma(t)$ belongs to the intersection $\mfd_A(M)\cap\mfd_A(N)$ of two walls, then the dimension vector $\dimv M$ of $M$ is a scalar multiple of the dimension vector $\dimv N$ of $N$;
\item  whenever $\gamma(t)$ is in $\mfd_A(M)$, then $\langle \gamma'(t),\dimv M\rangle\neq 0.$
\end{itemize}
 
That is, a smooth path is $\mfd_A$-generic if it crosses one wall at a time and the crossing is transversal.
\begin{lem}\label{l:mutations and generic paths}
    Let $(M,P)$ and $(N,Q)$ be basic $\tau$-tilting pairs. 
    \begin{itemize}
    \item[(1)] If $\gamma$ is a $\mfd_A$-generic path crossing finitely many walls such that $\gamma(0)\in \mfc_{(M,P)}$ and $\gamma(1)\in \mfc_{(N,Q)}$, then   $(M,P)\sim (N,Q)$. 
    \item[(2)] Conversely, if  $(M,P)\sim (N,Q)$, then there is a $\mfd_A$-generic path $\gamma$ crossing finitely many walls such that $\gamma(0)\in \mfc_{(M,P)}$ and $\gamma(1)\in \mfc_{(N,Q)}$.
    \end{itemize}
\end{lem}
\begin{proof}
Taking into account the description of chambers, the first statement is a direct consequence of Theorem \ref{l:mutation and neighboring  chamber}. The converse statement is proved in \cite{BST}.
\end{proof}

For the sake of distinction, for  a basic $\tau$-tilting pair $(M,P)$ of a finite-dimensional $k$-algebra $A$, we also denote the associated chamber by $\mfc_{(M,P)}^A$.
\begin{corollary}\label{c:induce-mutation-sequence-in-quotient}
    Let $A$ be a finite-dimensional $k$-algebra and $B$ a quotient algebra of $A$ with $|A|=|B|$. 
Let $(M_1,P_1),(M_2,P_2)\in \opname{\tau-tilt}A$ and $(N_1,Q_1),(N_2,Q_2)\in \opname{\tau-tilt}B$ such that $\mfc_{(M_1,P_1)}^A\subseteq \mfc_{(N_1,Q_1)}^B$ and $\mfc_{(M_2,P_2)}^A\subseteq \mfc_{(N_2,Q_2)}^B$. If  $(M_1,P_1)\sim (M_2,P_2)$, then  $(N_1,Q_1)\sim (N_2,Q_2)$.
\end{corollary}
\begin{proof}
    By Lemma \ref{l:B wall is A wall} and the definition of generic path, every $\mfd_A$-generic path $\gamma_A$ with $\gamma(0)\in\mfc_{(M_1,P_1)}^A$ and $\gamma(1)\in\mfc_{(M_2,P_2)}^A$ is also a $\mfd_B$-generic path. Now the statement follows from Lemma \ref{l:mutations and generic paths}.
\end{proof}

\begin{remark}
Let $\pi:A\twoheadrightarrow B$ be the canonical homomorphism in Corollary \ref{c:induce-mutation-sequence-in-quotient}. The induction functor $\pi_!=-\otimes_AB_B:\mod A\to \mod B$ induces a map $\pi_!:\opname{\tau-tilt}A\to \opname{\tau-tilt}B$ (cf. \cite[Corollary 2.4]{B20}).
By considering $g$-vectors, one can show that
\begin{itemize}
    \item $\mfc_{(M,P)}^A\subseteq \mfc_{\pi_!(M,P)}^B$ for any $(M,P)\in \opname{\tau-tilt}A$;
    \item If $(M,P)$ is a mutation of $(N,Q)\in \opname{\tau-tilt}A$, then either $\pi_!(M,P)=\pi_!(N,Q)$ or $\pi_!(M,P)$ is a mutation of $\pi_!(N,Q)$.
\end{itemize} 
This yields an alternative proof for Corollary \ref{c:induce-mutation-sequence-in-quotient}.

\end{remark}

\section{$\tau$-tilting graph of quasi-silted algebras}
In this section, we prove the main result of this note: the $\tau$-tilting graph of any quasi-silted algebra is connected. 

\subsection{Quasi-silted algebras}
Let $\mathcal{H}$ be a hereditary abelian category, and let $\der^b(\mathcal{H})$ denote the bounded derived category of $\mathcal{H}$ with suspension functor $[1]$. Recall that a rigid object $T \in \mathcal{H}$ is called \emph{tilting} if the condition $\Hom_{\mathcal{H}}(T, X) = 0 = \Ext_{\mathcal{H}}^1(T, X)$ implies $X = 0$. 

An object $T \in \der^b(\mathcal{H})$ is called \emph{silting} if $\Hom_{\der^b(\mathcal{H})}(T, T[>0]) = 0$ and $\opname{thick}(T) = \der^b(\mathcal{H})$. Furthermore, an object $T \in \der^b(\mathcal{H})$ is called a \emph{2-term silting complex} \cite{BZ16} over $\mathcal{H}$ if it satisfies the following two conditions:
\begin{itemize}
    \item[(S1)] $\Hom_{\der^b(\mathcal{H})}(T, M[i]) = 0$ for every $M \in \mathcal{H}$ and all $i \neq 0, 1$.
    \item[(S2)] $T$ is a silting object in $\der^b(\mathcal{H})$.
\end{itemize}

A finite-dimensional $k$-algebra $\Lambda$ is called \emph{quasi-tilted} if there exist a hereditary abelian category $\mathcal{H}$ and a basic tilting object $T \in \mathcal{H}$ such that $\Lambda \cong \End_{\mathcal{H}}(T)$. Analogously, $\Lambda$ is called \emph{quasi-silted} if there exist a hereditary abelian category $\mathcal{H}$ and a basic 2-term silting complex $T$ over $\mathcal{H}$ such that $\Lambda \cong \End_{\der^b(\mathcal{H})}(T)$. 

A finite-dimensional $k$-algebra $\Lambda$ is called \emph{shod} (small homological dimension) \cite{CL99} if every indecomposable $\Lambda$-module $M$ satisfies either $\pd_\Lambda M \leq 1$ or $\id_\Lambda M \leq 1$. A fundamental property of shod algebras is that their global dimension is at most $3$. In particular, a shod algebra $\Lambda$ is called \emph{strictly shod} if $\gld \Lambda = 3$.  

For an algebra $\Lambda$, we define the following subcategories of $\ind \Lambda$:
\begin{align*}
    \mathcal{L}_\Lambda &:= \{Y \in \ind \Lambda \mid \text{for every predecessor } X \text{ of } Y, \ \pd_\Lambda X \leq 1\}, \\
    \mathcal{R}_\Lambda &:= \{X \in \ind \Lambda \mid \text{for every successor } Y \text{ of } X, \ \id_\Lambda Y \leq 1\}.
\end{align*}

The following theorem summarizes several foundational results on shod algebras. Statement (1) was established in \cite{BZ16}, statements (2) and (3) were proven in \cite{CL99}, and statement (4) is due to \cite{HRS}.

\begin{theorem}\label{thm: shod and quasisilted}
Let $\Lambda$ be a connected finite-dimensional algebra.
\begin{itemize}
    \item[(1)] $\Lambda$ is shod if and only if it is quasi-silted;
    \item[(2)] $\Lambda$ is shod if and only if $(\add (\mathcal{R}_\Lambda \setminus \mathcal{L}_\Lambda), \add \mathcal{L}_\Lambda)$ is a split torsion pair in $\mod \Lambda$;
    \item[(3)] If $\Lambda$ is shod, then the ordinary quiver of $\Lambda$ has no oriented cycles;
    \item[(4)] $\Lambda$ is shod with $\gld \Lambda \leq 2$ if and only if it is quasi-tilted.
\end{itemize}
\end{theorem}

We will also need the following equivalent characterizations of strictly shod algebras; cf.\ \cite[Theorem II.1.14]{HRS} and \cite[Proposition 2.4]{CL99}.

\begin{lemma}\label{lem:strictly-shod}
Let $\Lambda$ be a finite-dimensional shod algebra. The following conditions are equivalent:
\begin{itemize}
    \item[(1)] $\Lambda$ is strictly shod;
    \item[(2)] $\mathcal{L}_\Lambda$ does not contain all indecomposable projective modules;
    \item[(3)] $\mathcal{R}_\Lambda$ does not contain all indecomposable injective modules.
\end{itemize}
\end{lemma}

The following result, established in \cite{DFP}, plays a central role in the proof of our main theorem.

\begin{theorem}{\cite[Theorem 3.9]{DFP}}\label{thm:tau-reduction-shod}
Let $\Lambda$ be a finite-dimensional shod algebra and $(U,R)$ a basic $\tau$-rigid pair of $\Lambda$. Then the $\tau$-reduction $\Lambda_{(U,R)}$ of $\Lambda$ with respect to $(U,R)$ is also a shod algebra. In particular, for any idempotent $e$ of $\Lambda$, the quotient algebra $\Lambda/\Lambda e\Lambda$ is a shod algebra.
\end{theorem}

 \subsection{Quasi-tilted algebras}\label{ss:quasi-tilted-algebras}
Let $\mh$ be a hereditary abelian category and $T$ a basic tilting object in $\mh$. Throughout this subsection, we denote by $C = \End_\mh(T)$ the quasi-tilted algebra associated with $T$. The following foundational properties are standard; see, for instance, \cite[Propositions~1.8, 1.9, and 1.10]{HR}.

\begin{lemma}\label{lem:well-known-facts}
\begin{itemize}
    \item[(1)] The tilting object $T$ induces a torsion pair $(\mt, \mf)$ in $\mh$, where $\mt = \Fac T$ and $\mf = \Sub \tau T$. Furthermore, there are fully faithful functors
    \[ F = \Hom_\mh(T,-) \colon \mt \to \mod C \quad \text{and} \quad G = \Ext^1_\mh(T,-) \colon \mf \to \mod C. \]
    \item[(2)] Let $\mathcal{X} = \Im G$ and $\mathcal{Y} = \Im F$. Then $(\mathcal{X}, \mathcal{Y})$ is a split torsion pair in $\mod C$; that is, every indecomposable $C$-module belongs to either $\mathcal{X}$ or $\mathcal{Y}$. 
    \item[(3)] For any $M \in \mathcal{Y}$, $\pd_C M \leq 1$. For any $N \in \mathcal{X}$, $\id_C N \leq 1$.
    \item[(4)] $\my$ is closed under $\tau_C$, and $\mx$ is closed under $\tau^{-1}_C$.
\end{itemize}
\end{lemma}

\begin{prop}\label{p:tau-reduction of tilted algebra}
Let $Z \in \my$ be an indecomposable  $\tau_C$-rigid $C$-module, and let $P({^\perp\tau_C Z})$ be the projective generator of ${^\perp\tau_C Z}$. Then $P({^\perp\tau_C Z}) \in \my$, and the endomorphism algebra $\End_C P({^\perp\tau_C Z})$ is quasi-tilted.
\end{prop}

\begin{proof}
Recall that the functor $F$ restricts to an equivalence $F:\mt\to \my$.
Let $X = F^{-1}(Z)$. Then $X$ is an indecomposable rigid object in $\mh$ lying in $\mt$. Consider the universal extension of $T$ by $X$: 
\begin{align}\label{gongshi: universal extension of X and T}
    0 \to T \to U \to X^t \to 0,
\end{align}
where the induced map $\Hom_{\mh}(X, X^t) \to \Ext_{\mh}^1(X, T)$ is an epimorphism. By \cite[Proposition~2.4]{HR}, $M = U \oplus X$ is a tilting object in $\mh$ that lies in $\mt$. Its image $F(M)$ is therefore a rigid $C$-module in $\my$ satisfying $\pd_C F(M) \leq 1$. Since $|F(M)| = |M| = |T| = |C|$, it follows that $F(M)$ is a tilting $C$-module. Noting that $Z = F(X)$, we see that $Z$ is a direct summand of $F(M)$ and that $F(M) \in {^\perp\tau_C Z}$.

Applying the functor $F = \Hom_{\mh}(T,-)$ to the short exact sequence (\ref{gongshi: universal extension of X and T}) yields the short exact sequence of $C$-modules:
\begin{align}\label{gongshi:universal extension of Z and C}
    0 \to C \to F(U) \to F(X^t) \to 0, 
\end{align}
where $F(X^t) \cong F(X)^t = Z^t$. For any $L \in {^\perp\tau_C Z}$, applying $\Hom_C(-, L)$ to (\ref{gongshi:universal extension of Z and C}) induces the exact sequence
\[ \Ext_C^1(Z^t, L) \to \Ext_C^1(F(U), L) \to \Ext_C^1(C, L) = 0. \]
Since $\pd_C Z \leq 1$, Auslander--Reiten duality yields $\Ext_C^1(Z^t, L) \cong \D\Hom_C(L, (\tau_C Z)^t) = 0$, which implies $\Ext_C^1(F(U), L) = 0$. Consequently, we obtain $\Ext_C^1(F(M), L) = 0$, proving that $F(M)$ is a projective object in ${^\perp\tau_C Z}$. By comparing ranks (cardinality), $F(M)$ is a projective generator of ${^\perp\tau_C Z}$. Let $M'$ be the basic tilting object obtained by decomposing $M$ and removing redundant direct summands, so that $M'$ contains $X$ as a summand. Then $F(M') = P({^\perp\tau_C Z})$, and 
\[ \End_C(P({^\perp\tau_C Z})) = \End_C(F(M')) \cong \End_{\mh}(M') \]
is a quasi-tilted algebra because $M' \in \mt$.
\end{proof}

\subsection{From quasi-tilted algebras to cluster-tilted algebras}
Keep the notations as in Section \ref{ss:quasi-tilted-algebras}.
Denote by $\tau$ the Auslander-Reiten translation of $\der^b(\mh)$. The cluster category of $\mh$ is defined as the orbit category $\mc:=\md^b(\mh)/\tauni[1]$, which admits a canonical triangle structure, i.e., the projection $\pi:\der^b(\mh)\to \mc$ is a triangle functor, cf. \cite{K05}.  We refer to \cite{BMRRT} for the definition related to cluster category.

The tilting object $T$ in $\mh$ naturally induces a cluster-tilting object in $\mc$. Denote by $B=\End_{\mc}(T)$  the  {\it cluster-tilted algebra} associated with $T$. According to \cite{ABS,Z06}, one has
\[
B\cong C\ltimes \Ext^2(\D(C),C).
\]
 Furthermore, there is a short exact sequence of $B$-modules:
\begin{align}
0\rightarrow \Ext^2(\D(C),C) \rightarrow  B \xrightarrow{\pi} C\rightarrow 0,
\end{align}
where the natural projection $\pi: B\rightarrow C$ is an algebra homomorphism and admits a section $\sigma: C\rightarrow B$ such that 
\begin{align}\label{gongshi: pi circ sigma is id}
    \pi\circ \sigma=1_C.
\end{align}
Along with the homomorphisms $\sigma$ and $\pi$, we have pairs of restriction functors and induction functors:
\[
\pi^\ast=-\otimes_{C}C_B:\mod C\to \mod B,\quad \pi_{!}=-\otimes_B C_C:\mod B\to \mod C,
\]
\[
\sigma^\ast=-\otimes_BB_C:\mod B\to \mod C,\quad \sigma_{!}=-\otimes_CB_B:\mod C\to \mod B.
\]
Furthermore, $\sigma^\ast\circ \pi^\ast=\mathbf{1}_{\mod C}$ and $\pi_{!}\circ \sigma_{!}=\mathbf{1}_{\mod C}$ since $\pi\circ \sigma=1_C$. It is easy to see that $\sigma_{!}$ preserves the projective modules. Indeed, let $f_1,\dots, f_n$ be a complete set of pairwise orthogonal primitive idempotents of $C$ and $B$. Denote by $P(i)=f_iC$ the indecomposable projective $C$-module associated with $f_i$. Then $\sigma_{!}(P(i))\cong f_iB$, the indecomposable projective $B$-module associated with $f_i$.

For notational simplicity, given any right $C$-module $M$, we identify $\pi^\ast(M)$ with $M$ when viewed as a right $B$-module.

\begin{lem}\cite[Proposition 4.2]{SS}\label{l: conditioin for M otimes C is M}
  For any $M \in \mod C$, we have $M \otimes_C B\cong M$ if and only if $\id_C M\leq 1$.
\end{lem}

\begin{lem}\cite[Propositions 3.2 and 3.3]{Zito}\label{l: tau-rigid between C and B}
    Let $M$ be a $\tau_C$-rigid module. 
    \begin{itemize}
        \item If $\id_C M\leq 1$, then $M$ is a $\tau_B$-rigid module.
        \item If $\pd_C \tau_C M\leq 1$, then $M\otimes_CB$ is a $\tau_B$-rigid module.
    \end{itemize}
\end{lem}

For a finite dimensional algebra $A$ and a $\tau$-rigid $A$-module $X$, the $g$-vector of $X$ in $\mod A$ is just the $g$-vector of the $\tau$-rigid pair $(X,0)$, which is simply written as $g^A_X$.
\begin{lem}\label{l: g-vector of B and C}
    Let $M$ be a $\tau_C$-rigid module and $M\otimes_CB$ be a $\tau_B$-rigid module, then $$g_M^C=g_{M\otimes_CB}^B.$$ 
\end{lem}
\begin{proof}
  Let  
  \begin{align}\label{g:projective resolusion of M as C module}
      P_1\rightarrow P_0\rightarrow M\rightarrow 0
  \end{align}
  be a minimal projective presentation of $M$ as $C$-module. Apply $-\otimes_CB$ to (\ref{g:projective resolusion of M as C module}),  we have an exact sequence 
  \begin{align}\label{g:projective resolusion of M as B module}
    P_1\otimes_CB\rightarrow P_0\otimes_CB\rightarrow M\otimes_CB\rightarrow 0.  
  \end{align}
 Since $\sigma_{!}$ preserves projectives,  (\ref{g:projective resolusion of M as B module}) is a projective presentation of $M\otimes_CB$ as $B$-module. Note that after applying $-\otimes_BC$ to (\ref{g:projective resolusion of M as B module}), we return to (\ref{g:projective resolusion of M as C module}) by $\pi_{!}\circ \sigma_{!}=\mathbf{1}_{\mod C}$. Hence, 
(\ref{g:projective resolusion of M as B module}) is a minimal projective presentation of $M\otimes_CB$ as $B$-module. Therefore, $g_M^C=g_{M\otimes_CB}^B.$
\end{proof}

\subsection{Main result}

We are now ready to state the main result of this paper.
 \begin{theorem}\label{thm:connected of quasi-silted algebra}
  The $\tau$-tilting graph of any  quasi-silted algebra is connected.
\end{theorem}
\begin{proof}
Let $A$ be a basic quasi-silted algebra. We prove the assertion by induction on $|A|$. The statement is trivial when $|A|=1$. Assume that $|A|=n>1$ and that the assertion holds for any basic quasi-silted algebra $\Lambda$ with $|\Lambda|<n$. 

We first establish the following claim:
\begin{quote}
    If $(M,P)$ and $(M',P')$ are basic $\tau$-tilting pairs of $A$ that share a nonzero direct summand, then $(M,P) \sim (M',P')$.
\end{quote}
Indeed, let $(U,R)$ be a common nonzero direct summand of $(M,P)$ and $(M',P')$. By Lemma~\ref{l:tau-reduction}, $E_{(U,R)}(M,P)$ and $E_{(U,R)}(M',P')$ are $\tau$-tilting pairs of the $\tau$-reduction $A_{(U,R)}$. Theorem~\ref{thm:tau-reduction-shod} implies that $A_{(U,R)}$ is a quasi-silted algebra. Since $|A_{(U,R)}| < n$, the induction hypothesis ensures that the $\tau$-tilting graph of $A_{(U,R)}$ is connected. Utilizing Lemma~\ref{l:tau-reduction} once more, we conclude that $(M,P) \sim (M',P')$.

To prove that the $\tau$-tilting graph of $A$ is connected, it suffices to show that $(M,P) \sim (A,0)$ for every basic support $\tau$-tilting pair $(M,P)$ of $A$. If $P \neq 0$, then $(M,P)$ and $(0,A)$ share a common direct summand (since $P$ and $A$ share an indecomposable projective summand), which yields $(M,P) \sim (0,A)$ by our claim. On the other hand, by Theorem~\ref{thm: shod and quasisilted}~(3), there exists an indecomposable projective $A$-module $Q$ such that $\Hom_A(Q, A/Q) = 0$. In particular, $(A/Q, Q)$ forms a basic $\tau$-tilting pair of $A$. Applying our claim again, we obtain $(A,0) \sim (A/Q, Q) \sim (0,A)$. Consequently, $(M,P) \sim (A,0)$ holds generally.

For the remainder of the proof, we may assume that $(M,0)$ is a basic $\tau$-tilting pair of $A$. We divide the argument into two cases depending on whether $A$ is strictly shod.

\noindent{\bf Case 1: $A$ is strictly shod}.
Recall that we have a split torsion pair $(\add (\mathcal{R}_A \setminus \mathcal{L}_A), \add \mathcal{L}_A)$ in $\mod A$. Hence, the $A$-module $M$ can be decomposed as $M \cong M_1 \oplus M_2$, where $M_1 \in \add (\mathcal{R}_A \setminus \mathcal{L}_A)$ and $M_2 \in \add \mathcal{L}_A$. We claim that $M_2 \neq 0$. Supposing otherwise, we would have $M \in \add (\mathcal{R}_A \setminus \mathcal{L}_A)$. Since $\mathcal{R}_A$ is closed under successors and $M$ is sincere, it follows that $\D(A)\in \add \mathcal{R}_A$, which directly contradicts Lemma~\ref{lem:strictly-shod}~(3). Thus, $M_2 \neq 0$. 

By Lemma~\ref{lem:strictly-shod}~(2), there exists an indecomposable projective $A$-module $Q$ satisfying $Q \in \add (\mathcal{R}_A \setminus \mathcal{L}_A)$. Because $(\add (\mathcal{R}_A \setminus \mathcal{L}_A), \add \mathcal{L}_A)$ is a torsion pair, we obtain $\Hom_A(Q, M_2) = 0$. As a consequence, $(M_2, Q)$ is a $\tau$-rigid pair of $A$. Let $(N, Q')$ be any basic  $\tau$-tilting pair of $A$ admitting $(M_2, Q)$ as a direct summand. From our initial claim and the previous discussion, we conclude that $(M,0) \sim (N, Q') \sim (A,0)$.

\noindent{\bf Case 2: $A$ is a quasi-tilted algebra}.
Let $\mh$ be a hereditary abelian category with a basic tilting object $T$ such that $A \cong \End_{\mh}(T)$. Let $B = A \ltimes \Ext^2_A(\D(A), A)$ be the associated cluster-tilted algebra. According to \cite{BMRRT, BKL, FG}, the $\tau$-tilting graph of $B$ is connected. Let $(\mx, \my)$ be the split torsion pair in $\mod A$ induced by $T$, and decompose $M$ as $M = M_\mx \oplus M_\my$, where $M_\mx \in \mx$ and $M_\my \in \my$. We analyze this setup via two subcases:

\noindent{\bf Subcase 2a:} $M_\my = 0$. This condition is equivalent to $M = M_\mx \in \mx$. By Lemma~\ref{lem:well-known-facts}, we have $\id_A M \leq 1$. Lemma~\ref{l: tau-rigid between C and B} ensures that $M$ is a $\tau_B$-rigid module. 
Since $|M| = n$, it follows that $M$ is a $\tau_B$-tilting module. By Lemmas~\ref{l: conditioin for M otimes C is M} and \ref{l: g-vector of B and C}, we conclude that the corresponding chambers coincide, i.e., $\mathfrak{C}_{(M,0)}^A = \mathfrak{C}_{(M,0)}^B$. Since the $\tau$-tilting graph of $B$ is connected, we have $(M,0) \sim (B,0)$ in $\mod B$. Finally, Corollary~\ref{c:induce-mutation-sequence-in-quotient} yields $(M,0) \sim (A,0)$.

\noindent{\bf Subcase 2b:} $M_\my \neq 0$. If $M_\my$ contains an indecomposable projective module as a direct summand, then $(M,0) \sim (A,0)$ follows immediately from our first claim. Now assume that $M_\my$ contains no nonzero projective direct summands, and let $Z$ be an indecomposable direct summand of $M_\my$. By Proposition~\ref{p:tau-reduction of tilted algebra}, $N := P({^\perp \tau_A Z}) \in \my$ is a $\tau_A$-tilting $A$-module. Since $(N,0)$ and $(M,0)$ share $Z$ as a common nonzero direct summand, our initial claim guarantees that $(M,0) \sim (N,0)$. It thus suffices to show that $(N,0) \sim (A,0)$.

Because $N \in \my$, we have $\tau_A N \in \my$. Applying Lemma~\ref{l: tau-rigid between C and B}, we see that $N \otimes_A B$ is a $\tau_B$-rigid module. Furthermore, following the proof of Lemma~\ref{l: g-vector of B and C}, we find $|N \otimes_A B| = |N|$. Hence, $N \otimes_A B$ is a $\tau_B$-tilting module. Let $N = N_1 \oplus \cdots \oplus N_n$ be its decomposition into indecomposable summands. Lemma~\ref{l: g-vector of B and C} yields $g_{N_i}^A = g_{N_i \otimes_A B}^B$ for each $1 \leq i \leq n$, which implies $\mathfrak{C}_{(N,0)}^A = \mathfrak{C}_{(N \otimes_A B, 0)}^B$. Since the $\tau$-tilting graph of $B$ is connected, we have $(N \otimes_A B, 0) \sim (B,0)$. Noting that $\mathfrak{C}_{(A,0)}^A = \mathfrak{C}_{(B,0)}^B$, we conclude from Corollary~\ref{c:induce-mutation-sequence-in-quotient} that $(N,0) \sim (A,0)$. This completes the proof.
\end{proof}

\begin{corollary}\label{cor:reachable-in-face}
Every quasi-silted algebra has the reachable-in-face property.
\end{corollary}

\begin{proof}
Let $A$ be a quasi-silted algebra. Suppose that $(M,P)$ and $(N,Q)$ are basic $\tau$-tilting pairs sharing a common nonzero direct summand $(L,R)$. We need to show that there exists a path connecting $(M,P)$ and $(N,Q)$ in the face $\mathcal{H}_{(L,R)}(\opname{\tau-tilt} A)$. 

Let $A_{(L,R)}$ be the $\tau$-reduction of $A$ with respect to $(L,R)$. By Theorem \ref{thm:tau-reduction-shod}, $A_{(L,R)}$ is also a quasi-silted algebra. Consequently, Theorem \ref{thm:connected of quasi-silted algebra} implies that the $\tau$-tilting graph $\mh(\opname{\tau-tilt} A_{(L,R)})$ is connected. On the other hand, by Lemma \ref{l:tau-reduction}, there is a graph isomorphism 
\[
\mh(\opname{\tau-tilt} A_{(L,R)}) \cong \mathcal{H}_{(L,R)}(\opname{\tau-tilt} A).
\]
It follows immediately that the face $\mathcal{H}_{(L,R)}(\opname{\tau-tilt} A)$ is connected, which yields the desired path between $(M,P)$ and $(N,Q)$.
\end{proof}

\section{$\tau$-tilting graph of $g$-tame algebras}
 In this section, we observe certain sufficient conditions under which the connectivity of $\tau$-tilting graphs is preserved by quotients of algebras. This yields several new classes of algebras whose $\tau$-tilting graphs are connected.

\subsection{Quotients}
Recall that any quotient of a basic finite-dimensional $k$-algebra remains basic. We begin with the following observation.
\begin{lem}\label{l:existence of e}
Let $A$ be a finite dimensional $k$-algebra and $B$ be a quotient algebra of $A$ with $|B|<|A|$. Then there exists a nonzero primitive idempotent element $e$ of $A$ such that $B$ is a quotient algebra of $A/\langle e\rangle$.
\end{lem}
\begin{proof}
   Suppose that $|A|=n$. Let $f: A \twoheadrightarrow B$ be the quotient homomorphism, and let ${e_1,...,e_n}$ be a complete set of primitive orthogonal idempotents for $A$. Let $f_i = f(e_i) \in B$ be the induced idempotents. If all $f_i \neq 0$, then ${f_1,...,f_n}$ would form a complete set of primitive idempotents for $B$, implying $|B|=n=|A|$, contradicting $|B|<|A|$. So there is some $i$ such that $f_i=0$, \ie $f(e_i)=0$. Let $\pi_i: A\rightarrow A/\langle e_i\rangle$ be the canonical epimorphism, there exists an epimorphism $g$ from $A/\langle e_i\rangle$ to $B$ such that $g\pi_i=f$. Hence $B$ is  a quotient algebra of $A/\langle e_i\rangle$.
\end{proof}

\begin{prop}\label{p:main-prop-2}
Let $A$ be a finite-dimensional $k$-algebra with connected $\tau$-tilting graph. For any primitive idempotent $e \in A$, the quotient algebra $A/\langle e\rangle$ has connected $\tau$-tilting graph.

 \end{prop}
\begin{proof}
     Consider the projective $A$-module $eA$ and note that $P(^\perp\tau (eA))=A$. Applying Lemma \ref{l:tau-reduction} to the basic $\tau$-rigid pair $(eA,0)$, we obtain an order-preserving bijection $$E_{(eA,0)}:\opname{\tau-tilt}_{(eA,0)}A\to \opname{\tau-tilt} A/\langle e\rangle.$$
Given any two $\tau$-tilting pairs $(M',P')$ and $(N',Q')$ of $A/\langle e\rangle$, let $$(M,P)=E_{(eA,0)}^{-1}(M',P')\text{ and } (N,Q)=E_{(eA,0)}^{-1}(N',Q').$$
Since the $\tau$-tilting graph of $A$ is connected, there is a path connecting $(M,P)$ and $(A,0)$ in $\mathcal{H}(\opname{\tau-tilt} A)$. As $(eA,0)$ is a common direct summand of $(M,P)$ and $(A,0)$, Lemma~\ref{l:cwz} implies that there is a path in $\mathcal{H}_{(eA,0)}(\opname{\tau-tilt}A)$ connecting $(M,P)$ and $(A,0)$. Similarly, there is a path in $\mathcal{H}_{(eA,0)}(\opname{\tau-tilt}A)$ connecting $(N,Q)$ and $(A,0)$. We conclude that there is a path connecting $(M',P')$ and $(N',Q')$ in $\mathcal{H}(\opname{\tau-tilt}A/\langle e\rangle)$ by Lemma \ref{l:tau-reduction}
\end{proof}

 \begin{prop}\label{p: B is connected if A is connected when A=B}
Let $A$ be a finite dimensional $k$-algebra and $B$ a quotient algebra of $A$ with $|B|=|A|$.  
Suppose that each $B$-chamber contains at least one $A$-chamber. If the  $\tau$-tilting graph of $A$ is connected, then the $\tau$-tilting graph of $B$ is also connected.

\end{prop}
\begin{proof}
The statement is a direct consequence of Corollary \ref{c:induce-mutation-sequence-in-quotient}.
\end{proof}


\subsection{$\tau$-tilting graphs of $g$-tame algebras}





Let $A$ be a  finite dimensional $k$-algebra with $|A|=n$. Recall from \cite[Definition 7.6]{AK} that $A$ is {\it $g$-tame} if $\overline{\mf(A)}=\R^n$, where \[\mf(A)=\bigcup_{(M,P)\in \opname{\tau-tilt}A}C_{(M,P)}.\]

We need the following property of $g$-tame algebras (cf. \cite[Proposition 3.11]{PYK} and \cite[Corollary 7.8]{AK}).
 \begin{lem}\label{l:g-tame under quotient}
Let $A$ be a finite dimensional $k$-algebra and $B$ a quotient algebra of $A$. If $A$ is $g$-tame, then so is $B$.
\end{lem}


We are now ready to state and prove the main result of this section.
 \begin{theorem}\label{t:connectness preserving under quotient}
Let $A$ be a $g$-tame algebra and $B$ a quotient algebra of $A$. 
 If the  $\tau$-tilting graph $\mathcal{H}(\opname{\tau-tilt}A)$ of $A$ is connected, then so is $\mathcal{H}(\opname{\tau-tilt}B)$. 
 \end{theorem}
\begin{proof}
If $|A|=|B|$, then every $B$-wall is an $A$-wall according to Lemma~\ref{l:B wall is A wall}. The $g$-tame condition ensures each $B$-chamber contains at least one $A$-chamber. Connectivity follows from Proposition \ref{p: B is connected if A is connected when A=B}.

So assume now that $|A|>|B|$. Lemma~\ref{l:existence of e} yields a nonzero idempotent $e\in A$ such that $B$ is a quotient of $A/\langle e\rangle$ with $|B|=|A/\langle e\rangle|$. 
Proposition \ref{p:main-prop-2} gives $\mathcal{H}(\opname{\tau-tilt} A/\langle e\rangle)$ is connected, and Lemma \ref{l:g-tame under quotient} shows $A/\langle e\rangle$ remains $g$-tame. The result then reduces to  the case $|A|=|B|$.
\end{proof}

The $\tau$-tilting graph  is known to be connected for various classes of algebras, including cluster-tilted algebras arising from hereditary abelian categories \cite{BMRRT,BKL,FG}, 2-Calabi-Yau tilted algebras originating from marked surfaces except closed surfaces with exactly one puncture \cite{QZ,Y}, gentle algebras \cite{FGLZ}, and skew-gentle algebras \cite{HZZ22}.
Note that $2$-Calabi-Yau tilted algebras arising from marked surfaces are tame algebras (cf. \cite{GLaS}) and it is clear that  skew-gentle algebras are tame algebras. According to \cite{PYK}, tame algebras are $g$-tame.
As a consequence of Theorem \ref{t:connectness preserving under quotient}, we obtain a large class of algebras with connected $\tau$-tilting graphs.



\begin{corollary}\label{c:connectness of some algebra}
Let $A$ be one of the following algebras:
\begin{enumerate}
\item A skew-gentle algebra;
\item A cluster-tilted algebra of tame type;
\item A 2-Calabi-Yau tilted algebra arising from a marked surface that is not closed with exactly one puncture.
\end{enumerate}
Then for any quotient algebra $B$ of $A$, the $\tau$-tilting graph $\mathcal{H}(\opname{\tau-tilt}  B)$ is connected.
\end{corollary}


The following corollary is immediate from Corollary~\ref{c:connectness of some algebra}.
\begin{corollary}\label{c:connectness of some algebra2}
 Let $A$ be a quotient algebra of a  hereditary algebra of  tame type.  Then
 the  $\tau$-tilting graph of $A$ is  connected.
 \end{corollary}

\bibliographystyle{plain}
\bibliography{ref}

\end{document}